\documentclass[11pt,epsf]{article}
\usepackage{srcltx}
\usepackage{graphicx}
\usepackage{theorem, amsmath, amssymb}
\newtheorem{Theorem}{Theorem}[section]
\newtheorem{Definition}[Theorem]{Definition}
\newtheorem{Proposition}[Theorem]{Proposition}

\newtheorem{Lemma}[Theorem]{Lemma}

\theoremstyle{remark}
\newtheorem{Remark}[Theorem]{Remark}

  {\begin{equation}\left\{\begin{aligned}}%
  {\end{aligned}\right.\end{equation}\ignorespacesafterend}%

\newenvironment{SDE*}%
  {\begin{equation*}\left\{\begin{aligned}}%
  {\end{aligned}\right.\end{equation*}\ignorespacesafterend}%

\usepackage{color}
\begin{document}
\title{
Asymptotic Perron's method 
and simple Markov strategies \\ 
in stochastic games and control
}

\author{Mihai S\^{\i}rbu \footnote{University of Texas at Austin,
    Department of Mathematics, 1 University Station C1200, Austin, TX,
    78712.  E-mail address: sirbu@math.utexas.edu. The research of
    this author was supported in part by the National Science
    Foundation under Grant    DMS 1211988. Any opinions, findings, and conclusions or recommendations expressed in this material are those of the authors and do not necessarily reflect the views of the National Science Foundation.}}
\maketitle
\begin{abstract} We introduce  a modification  of Perron's method, where semi-solutions are considered in a carefully defined asymptotic  sense. With this definition, we can  show, in a rather elementary way, that in a zero-sum game  or a control problem (with or without model uncertainty), the value function over all strategies coincides with the value function over  Markov strategies discretized in time. Therefore, there are always discretized Markov  $\varepsilon$-optimal strategies, 
(uniform with respect to the bounded initial condition). With a minor modification, the method produces a value and approximate saddle points for an asymmetric game of feedback strategies vs. counter-strategies. 
\end{abstract}


\section{Introduction}
The aim of the paper is to introduce the Asymptotic Perron's Method, i.e. constructing a solution of the Hamilton-Jacobi-Belman-Isaacs (HJBI) equation as the supremum/infimum of carefully defined {\bf asymptotic semi-solutions}. Using this method we show, in a rather elementary way, that the value functions of  zero-sum games/control problems can be (uniformly) approximated by some simple Markov strategies for the weaker player (the player in the exterior of the sup/inf or inf/sup). From this point of view, we can think of the method as an alternative to the shaken coefficients method of Krylov \cite{krylov-shaken-coeff} (in the case of only one player, under slightly different technical assumptions), or to the related method of regularization of solutions of HJBI's by \'Swi{\c{e}}ch in \cite{swiech-1} and \cite{swiech-2} (for control problems or games in Elliott-Kalton formulation). The method of shaken coefficients has been recently used to study games in Elliott-Kalton formulation in \cite{bouchard-nutz} under a convexity assumption (not needed here). 

While the result on zero-sum games (under our standing assumptions) is rather new, but certainly expected, the goal of the paper is  to present the method.
To the best of our knowledge, this modification of Perron's method does not appear in the literature. In addition, we believe that it  applies to more general  situations than we consider here, and either using a stochastic formulation (as in the present work) or an analytic one (see Remark \ref{analytic}). Compared to the method of shaken coefficients of Krylov, or to the regularization of solutions by \'Swi{\c{e}}ch, the analytic approximation of  the value function/solution of HJB by smooth approximate solutions  is replaced by the Perron construction. The careful definition of asymptotic semi-solutions  allows us to prove that such semi-solutions work well with Markov strategies.
The idea of restricting actions to a deterministic time grid is certainly not new, we just provide a method that works well for such strategies/controls.
 The arguments  display once again the robustness of the Perron construction, combined with viscosity comparison. There is basically a large amount of freedom in choosing the definition of semi-solutions, as long as they compare (as a simple consequence of their definition) to the  value functions. Here, we consider such {\bf asymptotic} semi-solutions.

 Perron's method and it's possible modifications seem to be rather useful in dynamic programming analysis. 
One could use  Ishii's arguments (see \cite{ishii})   to construct a viscosity solution and later upgrade it's regularity (if possible) for the purpose of verification (as in \cite{JS})  or employ a modification of Perron's method  over stochastic semi-solutions, if one does not expect smooth solutions (see \cite{bs-1}, \cite{bs-2} and \cite{bs-3}). Furthermore, with an appropriate definition of semi-solutions in the stochastic sense and a suitable model, one can even treat games and problems with model uncertainty, where dynamic programming is particularly cumbersome (see \cite{sirbu}, \cite{sirbu-3}).
The present work represents an additional  step  in this program 
 of using Perron's method and  some  possible offshoots to approach dynamic programming (DP) {\bf without having to first prove the dynamic programming principle} (DPP). This  step, unlike the previous (non-smooth) ones  actually provides a stronger conclusion than the usual Dynamic Programming (DP) approach.
More precisely, it proves the existence of approximate simple Markov strategies (like  \cite{krylov-shaken-coeff}, or \cite{swiech-1}, \cite{swiech-2}) something that is unclear from the DPP alone. 

  The formulation of the games/control problems we consider is resembling of   the formulation of deterministic games/control problems in the seminal work \cite{krasovskii-subbotin-88}. What we do here is to propose a novel method to study such models in stochastic framework.

Beyond the particular results obtained here, the (analytic version of the) method seems to also be useful in proving convergence of numerical schemes (but probably not the rate of convergence)
in the spirit  of Barles-Souganidis \cite{barles-souganidis}, or in the study of discretized games with mixed strategies as in \cite{buc-li-q}. We do not pursue this direction here.

Compared to the so called Stochastic Perron Method employed in \cite{sirbu} or \cite{bs-3}, the method we introduce here is quite different. The idea in \cite{sirbu} and \cite{bs-3} was to use {\bf exact} semi-solutions in  the stochastic sense then do Perron. In order to do so, the flexibility of stopping rules or stopping times was needed, leading to the possibility to complete the analysis {\bf only} over general feedback strategies (elementary for the purpose of well posed-ness of the state equation) or general predictable controls. Here, we propose instead to use  {\bf asymptotic} semi-solutions in the Perron construction. The flexibility on this end, allows us to work with a  deterministic time grid, resulting in approximation over Markov strategies. Obviously, the analytic part of the proof
(the ``bump-up/down'' argument)
is  similar  to \cite{sirbu} or \cite{bs-3}, but  those   parts of the proof were  already  similar to the (purely analytic) arguments of Ishii for viscosity Perron \cite{ishii}.


\section{Set-up and Main Results}
We use  the state equation  and the standing assumptions from \cite{sirbu}.

\vspace{.1in} 

\noindent {\bf The stochastic state system:}
\begin{equation}\label{eq:SDE}\left\{
\begin{array}{ll}
 dX_t=b(t,X_t,u_t, v_t)dt+\sigma (t, X_t,u_t, v_t)dW_t,   \\ 
 X_s=x \in \mathbb{R}^d.
 \end{array} \right.
 \end{equation}
 We assume that $u$ and $v$ belong to some compact metric spaces $(U, d_u)$ and $(V,d_V)$, respectively.
 For each $s$,  the problem comes with a fixed probability space $(\Omega, \mathcal{P}, \mathcal{F})$, a fixed filtration $\mathbb{F}=(\mathcal{F}_t)_{s\leq t\leq T}$ satisfying the usual conditions and a fixed Brownian  motion $(W_t)_{s\leq t\leq T}$, with respect to the filtration $\mathbb{F}$. We suppress the dependence on $s$ all over the paper. We emphasize that the filtration may be strictly larger than the saturated filtration generated by the Brownian motion. 
The coefficients $b:[0,T]\times \mathbb{R}^d\times U\times V\rightarrow \mathbb{R}^d$ and 
$\sigma :[0,T]\times \mathbb{R}^d\times U\times V\rightarrow \mathcal{M}^{d\times d'}$ satisfy the 

\pagebreak 

\noindent {\bf Standing assumptions:}
\begin{enumerate}
\item {\bf (C)}  $b, \sigma$  are jointly continuous on $[0,T]\times \mathbb{R}^d\times U\times V$

\item {\bf (L)} $b,\sigma$ satisfy a uniform local Lipschitz condition in $x$, i.e.
\begin{equation*}\label{Lip}
|b(t,x,u,v)-b(t,y,u,v)|+|\sigma (t,x,u,v)-\sigma (t,y,u,v)|\leq L(K) |x-y|\ \ 
\end{equation*}
$\forall\  |x|,|y|\leq K, t\in [0,T],\  u\in U, v\in V$
for some $L(K)<\infty$, and 
\item {\bf (GL)}
$b,\sigma$ satisfy a global linear growth condition in $x$
$$  
|b(t,x,u,v)|+|\sigma(t,x,u,v)|\leq C(1+|x|)$$
$\forall\  |x|,|y|\in \mathbb{R}^d, t\in [0,T],\  u\in U, v\in V$ for some $C<\infty.$
 \end{enumerate}
 Perron's method is a local method (for differential operators), so we only need local assumptions, except for the global growth which ensures non-explosion of the state equation and comparison for the Isaacs equation (see \cite{sirbu}).

Consider a   bounded and continuous reward (for the player $u$) function $g:\mathbb{R}^d\rightarrow \mathbb{R}$. Fix an initial time $s\in [0,T]$.  
We are interested in  the optimization problem
$$\sup_u\inf _v \mathbb{E}[g(X^{s,x;u,v}_T]$$
which is formally associated to  the lower Isaacs equation
\begin{equation}\label{isaacs}\left \{ \begin{array}{ll} -v_t-H^-(t,x,v_x, v_{xx})=0,\\
v(T,\cdot)=g(\cdot).\end{array}\right.\end{equation}
Above, we use the notations
$$H^-(t,x,p,M)\triangleq \sup_{u\in U}\inf _{v\in V} L(t,x,u,v,p,M),$$
$$L(t,x,u,v,p,M)\triangleq b(t,x,u,v)\cdot p+\frac 12 Tr\big (\sigma (t,x,u,v)\sigma ^T(t,x,u,v)M\big).$$
It is well known that making sense of what $u$, $v$ should be above is a non-trivial issue, aside from relating the optimization problem to the Isaacs equation. We expect that three possible models are actually represented as  solutions of the lower  Isaacs equation. They are
\begin{enumerate}
\item the lower value of a symmetric game over feedback strategies (as in  \cite{pz-game} or \cite{sirbu}),
\item the value function of a robust control problem where $u$ is an intelligent maximizer and $v$ is a possible worst case scenario modeling Knigthian uncertainty (see \cite{sirbu-3}), or
\item the genuine value of a sup-inf/inf-sup non symmetric game over feedback strategies vs. feedback counter-strategies (as  in  \cite{krasovskii-subbotin-88}[Section 10] or \cite{fleming-hernandez2}, \cite{fleming-hernandez2-2}).
\end{enumerate}
Although the main goal of the paper is to present the novel modification of Perron's method, we also want provide a unified treatment for the three models above. Therefore, some modeling considerations need to be taken into account (on top of the ones we simply  borrow from \cite{sirbu} or \cite{sirbu-3}). 
More precisely, depending on the information structure, one has two types of zero-sum games. The first, and fully symmetric one, is where both players observe the past of the state process $X$ (only), and therefore, make decisions based on this. In this model, both players use feedback strategies, i.e.  non-anticipative 
functionals of the path $X_{\cdot}$ (this is the case in \cite{pz-game} or \cite{sirbu}). On the other hand, one can envision a game where player $u$ can only see the state, but the player $v$ observes the state, and, in addition, the control $u$ (in real time).  The intuition corresponding to the nature of the problem, tells us that the advantage that the player $v$ can gain from observing the whole past of the control $u$ actually  comes {\bf only }  from being able to adjust, instantaneously to the observation $u_t$. In other words, in such a model we use  a (counter) strategy for player $v$ that depends on
\begin{enumerate}
\item the whole past of the state process $X$ up to the present time $t$
\item (only)  the current value of the adverse control $u_t$.
\end{enumerate}
This modeling avenue for deterministic games is taken in Section 10 of  the seminal monograph \cite{krasovskii-subbotin-88} and followed up in the important work on stochastic games \cite{fleming-hernandez2}, \cite{fleming-hernandez2-2}.
\begin{Definition}[Feedback Strategies and Counter-Strategies]  Fix a time $s$.

\begin{enumerate}
\item  a feedback strategy for player $u$ is a mapping 
$$\alpha :[s,T]\times C[s,T]\rightarrow  U,$$
which is  predictable with respect to the (raw) filtration on the path space
$$\mathcal{B}_t\triangleq\sigma (y(q),s\leq q\leq t),\ \ \forall t\in [s,T].$$
A similar definition holds for the player $v$.
\item a feedback counter-strategy for the player $v$ is a mapping
$$  \gamma:[s,T]\times C[s,T]\times U \rightarrow V,$$ 
which is measurable with respect to $\mathcal{P}^s\otimes \mathcal{U}^b/\mathcal{V}$. The second player uses, at time $t$, the action
$$v_t=\gamma (t, X_{\cdot}, u_t),$$
where $u_t$ is the action of the player $u$ at time $t$. Above, 
$\mathcal{P}^s$ is the predictable sigma-field on $[s,T]\times C[s,T]$ with respect to the raw filtration on the path space and $\mathcal{U}^b$ is the Borel sigma-field on $U$ with respect to the metric $d_U$.

\end{enumerate}
\end{Definition}
The definition of counter-strategies goes back to  \cite{krasovskii-subbotin-88} (see, for example, the Definition on page 431 for a Markovian version)  and the name is also borrowed from there. Feedback strategies have been used in control/games for a long time, and it is hard to trace their exact  origin.
It is clear that, for a fixed  pair of  feedback strategies $(\alpha, \beta)$ (or a feedback strategy vs. a feedback counter-strategy $(\alpha, \gamma)$) the state equation may fail to have a solution. Therefore, we need some restrictions on strategies (and counter-strategies). We use here (for strategies,  there is) the restriction to elementary strategies in \cite{sirbu}.

\begin{Definition}[Elementary feedback strategies and counter-strategies] Fix $s$.
\begin{enumerate}
\item (see \cite{sirbu}) an elementary feed-back strategies (for the $u$-player) is a  predictable functional  of the path
$\alpha :[s,T]\times C[s,T]\rightarrow  U ,$
which is only "rebalanced" at some  stopping rules $\tau _1\leq \tau _2\leq ...\leq \tau _n$. More precisely, there exists $n$ and  some
$$\tau _k:C[s,T]\rightarrow [s,T], \ \ k=1,2, \dots, n \ \ \textrm{such \ that}\ \   \forall k, \ \ \ \{\tau _k\leq t\}\in \mathcal{B}_t, \ s\leq t\leq T,$$
and there also exist some  $$\xi_k :C[s,T]\rightarrow  U,\ \ \textrm{with }
\ \ \xi _k  \in \mathcal{B}_{\tau_{k-1}}$$ so that 
$$\alpha (t,y)=\sum _{k=1}^n 1_{\{\tau _{k-1}(y)<t\leq \tau _k(y)\}} \xi _k(y), \ \forall \ s\leq t\leq T, \ y\in C[s,T].$$
A similar definition holds for player $v$.
 We denote by $\mathcal{A}(s)$ and  $\mathcal{B}(s)$ the collections of {\bf elementary (pure) feed-back} strategies for the $u$-player and the $v$-player, respectively.

\item  an elementary feed-back counter-strategy for the $v$-player is a 
mapping
$$  \gamma:[s,T]\times C[s,T]\times U \rightarrow V,$$ for which there exist 
stopping rules $\tau _1\leq \tau _2\leq ...\leq \tau _n$ as above  and some  $$\eta_k :C[s,T]\times U \rightarrow  V,\ \ \textrm{with }
\ \ \eta _k  \in \mathcal{B}_{\tau_{k-1} } \otimes \mathcal{U}^b$$ 
 so that 
$$\gamma (t,y,u)=\sum _{k=1}^n 1_{\{\tau _{k-1}(y)<t\leq \tau _k(y)\}} \eta _k(y,u),\ \forall \ s\leq t\leq T, \ y\in C[s,T],\ u\in U.$$
We denote by $\mathcal{C}(s)$ the set of elementary counter-strategies of the $v$-player.

 \end{enumerate}
\end{Definition}
We also denote by $\mathcal{U}(s)$ and  $\mathcal{V}(s)$ the set of {\bf open-loop controls} for the $u$-player and the $v$-player, respectively. Precisely, 
 $$\mathcal{V}(s)\triangleq \{v:[s,T]\times \Omega \rightarrow V|\ 
\textrm{predictable \ with\  respect\  to\  } \mathbb{F}=(\mathcal{F}_t)_{s\leq t\leq T}\},$$
and a similar definition is made for  $\mathcal{U}(s)$. Since the number of symbols may become overwhelming, we will use for notation
\begin{enumerate}
\item  $\alpha, \beta$ for the feedback strategies of players $u$ and $v$,
\item  $u, v$ for the open loop controls,
\item $\gamma$ for the feedback counter-strategy of the second player $v$.
\end{enumerate}
The {\bf only} reason  to restrict feedback strategies or counter-strategies to be elementary, is to have well posed-ness of the state equation (actually in the strong sense).
\begin{Lemma}\label{lemma:well-posed} Fix $s,x$. Assume that the first player uses an open-loop control $u\in \mathcal{U}(s)$ or an elementary feedback strategy $\alpha \in \mathcal{A}(s)$. Assume also that  the second player uses either an open-loop control $v \in \mathcal{V}(s)$ or a feedback strategy $ \beta \in  \mathcal{B}(s)$ or a counter-strategy $\gamma \in \mathcal{C}(s)$ . Then the state equation has a unique strong solution.
\end{Lemma} The result has been briefly proven in \cite{sirbu} and \cite{sirbu-3} with the possible exception of both players using open-loop controls (where the result is both obvious, and not really used here) and the case when the second player uses a counter strategy $\gamma$. We treat here the case when the second player uses an elementary counter-strategy $\gamma$.
\begin{enumerate}
\item $(u, \gamma)$. If an open-loop control $u$ and a counter-strategy $\gamma$ are fixed, one has to solve the state equation, iteratively, in between stopping rules $\tau _{k-1}, \tau_{k}$. This is possible since, in between this stopping rules, the state equation has no ``feedback'', but it really looks like an open-loop vs. open-loop equation. In other words,  assume that the state equation is well posed up to the stopping rule $\tau _{k-1}$ and the value of the process at that time  (i.e. the random variable  $X_{\tau _{k-1}(X_{\cdot})}$) is $L^2$). Once the state process "has arrived" to the stopping rule $\tau _{k-1}$ one can start solving  the SDE with an initial random time  (rather than stopping rule) $\tau _{k-1}(X_{\cdot})$  and initial condition $X_{\tau _{k-1}(X_{\cdot})}$ and using the predictable open-loop controls
$$(u_t)_{\tau _{k-1}(X_{\cdot}) \leq t \leq T},\ \ \ (\eta _k ( X_{\cdot}  ,u_t))_{\tau _{k-1}(X_{\cdot}) \leq t \leq T},$$
on the stochastic interval $\tau _{k-1}(X_{\cdot})\leq t\leq T$. 
 This can be solved due to  the locally Lip and linear growth assumptions and results in a square integrable process $X$ defined from $\tau _{k-1}(X_{\cdot})$ forward. We let this run up to  the next stopping rule $\tau _k$ and then continue as above. The result is a unique strong and square integrable solution $X$.
\item $(\alpha, \gamma)$ If an elementary feedback strategy $\alpha \in \mathcal{A}(s)$ and an elementary  counter-strategy $\gamma \in \mathcal{C}(s)$ are chosen, one has to simply see that the super-position of the counter-strategy over the strategy is an elementary strategy. In other words
$$\gamma  [\alpha ] (t, y)\triangleq \gamma (t, y, \alpha (t, y))\ \ \forall s\leq t\leq T,\  y\in \mathcal{C}[s,T]$$
defines a $\gamma  [\alpha]\in \mathcal{B}(s)$. For $\alpha, \gamma [\alpha]$ the state equation is well-posed (see \cite{sirbu}) therefore it is well posed over $\alpha, \gamma$.
\end{enumerate}
  The state process $X$ will be used with some explicit (and obvious) super-script notation $X^{s,x;\cdot, \cdot}$.
Following \cite{sirbu} we first define the lower value function (only) for the symmetric game in \cite{sirbu} in between two feedback players:\begin{equation}\label{sym-game}  V^-(s,x)\triangleq  \sup _{\alpha \in \mathcal{A}(s)}\left (\inf _{\beta \in \mathcal{B}(s)}  \mathbb{E}[g(X^{s,x;\alpha ,\beta}_T)]\right). 
\end{equation} 
Next, following \cite{sirbu-3} we consider  a robust control problem where the intelligent player $u$ uses feedback strategies and the  open-loop controls $v$ parametrize worst case scenarios/Knightian uncertainty:
\begin{equation}\label{robust-control} v^-(s,x)\triangleq  \sup _{\alpha \in \mathcal{A}(s)}\left (\inf _{v \in \mathcal{V}(s)}  \mathbb{E}[g(X^{s,x;\alpha ,v}_T)]\right).
\end{equation} 
We recall that the filtration may be larger than the one generated by the Brownian motion. 
Finally,  following \cite{krasovskii-subbotin-88} or \cite{fleming-hernandez2}  we consider a genuine game (with a lower and an upper value) in between two intelligent players who can both observe the state, but the second player has the advantage of also observing first player's actions in real time, i.e. 
\begin{equation}\label{non-sym-game}
W^-(s,x)\triangleq  \sup _{\alpha \in \mathcal{A}(s)}\left (\inf _{\gamma \in \mathcal{C}(s)}  \mathbb{E}[g(X^{s,x;\alpha ,\gamma }_T)]\right) \leq 
   \inf _{\gamma  \in \mathcal{C}(s)}\left (\sup _{\alpha \in \mathcal{A}(s)}  \mathbb{E}[g(X^{s,x;\alpha ,\gamma }_T)]\right) \triangleq W^+(s,x).
\end{equation}
In addition to this, for mathematical reasons, we define yet another value function
\begin{equation}\label{v+}
v^+(s,x)\triangleq \inf _{\gamma  \in \mathcal{C}(s)}\left (\sup _{u \in \mathcal{U}(s)}  \mathbb{E}[g(X^{s,x;u ,\gamma }_T)]\right )\geq W^+(s,x).
\end{equation}
We could attach to $v^+$ the meaning of some robust optimization   problem, but this is not natural, since the intelligent optimizer $v$ can see in real time the ``worst case scenario''.
By  simple observation we  have  
$v^-\leq W^-\leq W^+ \leq v^+.$ In addition, since $\mathcal{B}(s)\subset \mathcal{C}(s)$, and according the to second part of the proof of Lemma \ref{lemma:well-posed}, for any fixed $\alpha \in \mathcal{A}(s), \gamma \in \mathcal{C}(s)$ we have $\gamma [\alpha]\in \mathcal{B}(s)$, we actually see  that 
$$     \sup _{\alpha \in \mathcal{A}(s)}\left (\inf _{\gamma \in \mathcal{C}(s)}  \mathbb{E}[g(X^{s,x;\alpha ,\gamma }_T)]\right)
\leq  \sup _{\alpha \in \mathcal{A}(s)}\left (\inf _{\beta \in \mathcal{B}(s)}  \mathbb{E}[g(X^{s,x;\alpha ,\beta}_T)]\right)\leq  \sup _{\alpha \in \mathcal{A}(s)}\left (\inf _{\gamma \in \mathcal{C}(s)}  \mathbb{E}[g(X^{s,x;\alpha ,\gamma [\alpha] }_T)]\right),$$
i.e. $W^-=V^-$. Altogether,  we know that
$$v^-\leq W^-=V^-\leq W^+\leq v^+.$$
\begin{Remark} Since the state equation is well posed, the counter-strategies $\gamma\in \mathcal{C}(s)$ (and, therefore, the feedback strategies $\beta \in \mathcal{B}(s)$) are all strategies in the sense of Elliott-Kalton \cite{ek} (or \cite{fs} for the stochastic case). Therefore, there is a natural question: since, in the non-symmetric game \eqref{non-sym-game} the player $v$ observes $u$, why not formulate the game as a sup/inf and inf/sup over Elliott Kalton strategies vs. open-loop controls?
In other words, we could set up the problem (with little rigor in   formulation) as
$$\sup _{u }\inf _{e}\mathbb{E}[g(X^{s,x;u,e}_T)]\leq \inf _{e} \sup _{u }\mathbb{E}[g(X^{s,x;u,e}_T)]
$$
where $e$ is an Elliott-Kalton strategy. The lower value of the game above is, heuristically (quite obvious in the deterministic case), equal to the lower value of the symmetric game over open-loop controls and the upper value above is expected (according to  \cite{fs}) to be the solution of the lower  Isaacs equation (or, the unified Isaacs equation if the Isaacs condition holds). Therefore, (a possible modification of) the well known example of Buckdahn (Example 8.1 in \cite{pz-game}) shows that such game may fail to  have a value.
 On the other hand, the non-symmetric game over feedback strategies vs. counter-strategies will have a value (see Theorem \ref{main} below), showing that the information structure considered in \eqref{non-sym-game} is better suited to analyze a non-symmetric game of this type.  Elliott-Kalton strategies are designed to be considered only in the exterior of the inf-sup/sup-inf, i.e. only in the upper value above. 

\end{Remark}

We now define a  very special class of elementary strategies, namely, simple Markov strategies. Actually, what we call simple Markov strategies below are called "positional strategies " in \cite{krasovskii-subbotin-88} (see the Definition on page 12 and relation (3) on page 6, for example) and are used extensively in deterministic games. Similar  strategies/counter-strategies  are used under different names in the more recent interesting contributions \cite{fleming-hernandez2}, \cite{fleming-hernandez2-2} on asymmetric zero-sum games.

\begin{Definition}[time grids, simple  Markov strategies and counter-strategies]  Fix $0\leq s\leq T$.
\begin{enumerate}

\item A time grid for $[s,T]$ is a finite sequence $\pi$ of $s=t_0<t_1<\dots<t_n=T$.

\item Fix a time grid $\pi$ as above. A strategy $\alpha \in \mathcal{A}(s)$ is called  a simple Markov strategy {\bf over} $\pi$  if there exist
 some functions $\xi_k:\mathbb{R}^d\rightarrow U, k=1,\dots, n$ measurable, such that
$$\alpha (t, y(\cdot))=\sum _{k=1}^n 1_{\{t_{k-1}<t\leq t_k\}} \xi _k(y(t_{k-1})).$$
The set of all simple  Markov strategies over  $\pi$ is denoted by $\mathcal{A}^M(s,\pi).$
Define the set of all simple Markov strategies, over all possible time grids as 
$$\mathcal{A}^M(s)\triangleq \bigcup _{\pi} \mathcal{A}^M(s,\pi).$$

\item Fix a time grid $\pi$ as above. A counter-strategy $\gamma \in \mathcal{C}(s)$ is called  a simple Markov counter-strategy {\bf over} $\pi$  if there exist
 some functions $\eta_k:\mathbb{R}^d\times U\rightarrow V, k=1,\dots, n$ measurable, such that
$$\gamma (t, y(\cdot), u)=\sum _{k=1}^n 1_{\{t_{k-1}<t\leq t_k\}} \eta  _k(y(t_{k-1}), u).$$
The set of all simple  Markov counter strategies over  $\pi$ is denoted by $\mathcal{C}^M(s,\pi).$
Define the set of all simple Markov  counter-strategies, over all possible time grids as $$\mathcal{C}^M(s)\triangleq \bigcup _{\pi} \mathcal{C}^M(s,\pi).$$
\end{enumerate}

\end{Definition}
In words,  for a strategy "simple Markov" means that the player only changes actions over the time grid, and anytime he/she does so, the new control depends on the current position only. However, for ``simple  Markov counter-strategies'', the situation stands in stark contrast: actions are changed continuously in time based on the instantaneous observation of $u$, but the information from observing the state is only updated discretely over $\pi$.
\begin{Remark}\label{counter-discrete} If we attempt to fully discretize counter-strategies as well, allowing for the $v$-player to only update his/her actions  over the time grid $\pi$, one cannot expect the game \eqref{non-sym-game} to have a value, since the opponent $u$ can change actions many times in between $t_{k-1}$ and $t_k$. One would need to restrict both players to the same time grid $\pi$ then to pass to the limit over $\pi$ (in the spirit of \cite{buc-li-q}) to obtain a value over fully discretized counter-strategies.
\end{Remark}
Consider now the same optimization problems  as above, but where  the weaker player is restricted to using {\bf only} simple Markov strategies/counter-strategies, restricted to a fixed time grid, or not. Denote by 
$$ v _{\pi} ^-(s,x)\triangleq  \sup _{\alpha \in \mathcal{A}^M(s,\pi)}\left (\inf _{v \in \mathcal{V}(s)}  \mathbb{E}[g(X^{s,x;\alpha ,v}_T)]\right) \leq v^-(s,x)\leq  W^-(s,x)= V^-(s,x),$$
and 
$$ v^- _M(s,x)\triangleq  \sup _{\alpha \in \mathcal{A}^M(s)}\left (\inf _{v \in \mathcal{V}(s)}  \mathbb{E}[g(X^{s,x;\alpha ,v}_T)]\right) =\sup _{\pi} v_{\pi}^- (s,x) \leq v^-(s,x)\leq W^-(s,x)= V^-(s,x),$$
as well as 
$$ v^+ _{\pi}(s,x)\triangleq  \inf _{\gamma \in \mathcal{C}^M(s,\pi)}\left (\sup  _{u \in \mathcal{U}(s)}  \mathbb{E}[g(X^{s,x;u ,\gamma}_T)]\right) \geq v^+(s,x)\geq W^+(s,x),$$
and 
$$ v^+ _M(s,x)\triangleq  \inf _{\gamma  \in \mathcal{C}^M(s)}\left (\sup _{u \in \mathcal{U}(s)}  \mathbb{E}[g(X^{s,x;u ,\gamma}_T)]\right) =\inf _{\pi} v_{\pi}^+ (s,x) \geq v^+(s,x)\geq W^+(s,x).$$

The main result is that the $u$-player (in either setting) cannot do better with general feedback strategies than with simple Markov strategies, and the $v$ player can do as well over simple Markov counter-strategies as over general (elementary) ones.  The non-symmetric game \eqref{non-sym-game} has a value and admits approximate saddle points over simple  Markov strategies/simple Markov counter-strategies. This is certainly expected. Under slightly stronger technical assumption (natural filtration, global Lipschitz conditions) this is proven in \cite{fleming-hernandez2} with different methods, and is also pursued in Section 10 of \cite{krasovskii-subbotin-88} for deterministic games. The one player (i.e. control) case is well studied (again with different methods) in the seminal monograph \cite{MR2723141} (see subsection \ref{one-player} for more comments). 

The  whole idea of approximating the problem over controls/strategies that are fixed over time grids is by no means novel. It  goes much further than   \cite{krasovskii-subbotin-88} for games or  \cite{MR2723141} in control. The main contribution of the paper is to obtain such results with a different  (and more elementary) method (and  under slightly different assumptions).
\begin{Theorem}\label{main} Under the standing assumptions, we have that
$$v_M^-=v^-=W^-=V^-=W^+ =v^+ =v^+_M$$ and the common value is the unique bounded  continuous viscosity solution of the lower Isaacs equation. In particular, the non-symmetric game \eqref{non-sym-game} has a value which is equal to the lower value of the symmetric game \eqref{sym-game}. For each $N$, and each $\varepsilon>0$, there exists $\delta (N,\varepsilon)>0$ such that
$$\forall s\in [0,T], \forall \ |\pi|\leq \delta,\ \exists \ \hat{\alpha}  \in \mathcal{A}^M(s,\pi), \hat{\gamma } \in \mathcal{C}^M(s,\pi)$$
such that $\forall\ |x|\leq N,$ we have 
$$0\leq  W^-(s,x)-\underbrace{\inf _{v \in \mathcal{V}(s)}  \mathbb{E}[g(X^{s,x;\hat{\alpha} ,v}_T)]}_{\leq v^-_{\pi}(s,x)}\leq \varepsilon  \ \textrm{and}\ 
0\leq \underbrace{\sup _{u\in \mathcal{U}(s)}\mathbb{E}[g(X^{s,x;u,\hat{\gamma}}_T)]}_{\geq v^+_{\pi}(s,x)} -W^+(s,x) \leq \
\varepsilon.$$
Therefore, 
$v^-_{\pi}, v^+_{\pi}\rightarrow V^-=W^-=W^+$ as $|\pi|\rightarrow 0$ uniformly on compacts in $[0,T]\times \mathbb{R}^d$.
\end{Theorem}
 The result above can be rewritten as 
 $$\mathbb{E}[g(X^{s,x;u, \hat{\gamma}}_T]-\varepsilon\leq W^+(s,x)=V^-(s,x)=W^-(s,x) \leq 
 \mathbb{E}[g(X^{s,x;\hat{\alpha} , v }_T]
  +\varepsilon\  (\forall) \ \ (u,v) \in \mathcal{U}(s)\times \mathcal{V}(s), |x|\leq N,$$
(this is the way the result is phrased in the very interesting paper \cite{fleming-hernandez2}). Plugging in above the  open-loop controls 
$$\hat{u}_t=\hat{\alpha}(t, X^{s,x;\hat {\alpha}, \hat {\gamma}}_{\cdot}),  \hat{v}_t=\hat{\gamma}(t, X^{s,x;\hat {\alpha}, \hat {\gamma}}_{\cdot}, \hat{u}_t)$$
we obtain that 
$|\mathbb{E}[g(X^{s,x;\hat{\alpha}, \hat{ \gamma}}_T]-V^-(s,x)|\leq \varepsilon.$ This yields
 $$\mathbb{E}[g(X^{s,x;u, \hat{\gamma}}_T]-2\varepsilon\leq \mathbb{E}[g(X^{s,x;\hat{\alpha}, \hat{ \gamma}}_T]
 \leq 
 \mathbb{E}[g(X^{s,x;\hat{\alpha} , v }_T]
  +2 \varepsilon\  (\forall) \ \ (u,v) \in \mathcal{U}(s)\times \mathcal{V}(s), |x|\leq N,$$
  which is stronger than 
  $$\mathbb{E}[g(X^{s,x;\alpha , \hat{\gamma}}_T]-2\varepsilon\leq \mathbb{E}[g(X^{s,x;\hat{\alpha}, \hat{ \gamma}}_T]
 \leq 
 \mathbb{E}[g(X^{s,x;\hat{\alpha} , \gamma }_T]
  +2 \varepsilon\  (\forall) \ \ (\alpha ,\gamma) \in \mathcal{A}(s)\times \mathcal{C}(s), |x|\leq N.$$
This means that  $(\hat{\alpha}, \hat{\gamma})$ is actually  a $2\varepsilon$-saddle point for the (genuine) non-symmetric game \eqref{non-sym-game}. 

\section{Proof: the Asymptotic Perron's Method}
We introduce here the new version of Perron's method, where semi-solutions of the (lower) Isaacs equations are replaced by {\bf asymptotic} semi-solutions.  In our particular framework, we use asymptotic {\bf stochastic} semi-solutions (so the Perron Method here is asymptotic in the stochastic sense of \cite{sirbu} or 
\cite{bs-3}). However, we claim that a similar Asymptotic Perron Method can be designed in the {\bf analytic} framework (see Remark \ref{analytic}). 

The definition of asymptotic semi-solutions is  different from the definition of stochastic semi-solutions in \cite{sirbu-3} or \cite{sirbu}, and, consequently, so are the proofs. 
The analytic part of the proof still resembles Ishii \cite{ishii} and the probabilistic part uses It\^o along the smooth test function, but this is where similarities stop.
 As mentioned, the  method we introduce here, since it comes in close relation to Markov strategies, can be viewed as an alternative to the powerful method of shaken coefficients of Krylov, 
\cite{krylov-shaken-coeff} or  to the work of \'Swi{\c{e}}ch  \cite{swiech-1} and \cite{swiech-2}. Anyway, in the case  of games, the method of  \'Swi{\c{e}}ch actually works for a slightly different game, defined on a space accommodating an independent Brownian motion, and using Elliott-Kalton strategies. Since we are ultimately studying a non-symmetric game \eqref{non-sym-game}, the analysis has to be done separately for the two value functions. 

\subsection{Asymptotic Perron over strategies}\label{subsection-strategies}

We perform here a sup Perron construction lying below the value function $v_M^-$. If one only cares about the lower value of the game $V^-$ in \eqref{sym-game} or the robust control problem \eqref{robust-control}, this is the only construction we need. Together with the results in \cite{sirbu} this provides a full approximation of the two problems by elementary Markov strategies of the player $u$. We regard   this as the most important result.

\begin{Definition}[Asymptotic Stochastic Sub-Solutions]\label{def-sub-sol}
A function $w:[0,T]\times \mathbb{R}^d\rightarrow \mathbb{R}$ is called an asymptotic (stochastic) sub-solution of the (lower) Isaacs equation, if it is bounded, continuous and  satisfies $w(T,\cdot)\leq g(\cdot)$. In addition, there exists a gauge function $\varphi=\varphi_w:(0,\infty)\rightarrow (0,\infty),$ depending on $w$ such that
\begin{enumerate}
\item $\lim _{\varepsilon\searrow 0}\varphi(\varepsilon)=0,$
\item for each $s$ (and the optimization problem coming with it),  for each time $s\leq r\leq T$, there exists a measurable function $\xi :\mathbb{R}^d\rightarrow U$ such that, for each $x$, each  $\alpha \in \mathcal{A}(s)$ and $v\in \mathcal{V}(s)$, if we make the notation
$\alpha [r,\xi]\in \mathcal{A}(s),$
defined by
$$\alpha [r,\xi](t, y(\cdot))=1_{\{s<t\leq r\}}\alpha (t, y(\cdot))+1_{\{r<t\leq T\}} \xi(y(r)),$$
then, for each $r\leq t\leq T$ we have
\begin{equation}\label{sub-sol}
w(r, X^{s,x; \alpha,v}_r)=w(r, X^{s,x; \alpha[r,\xi],v}_r)\leq 
\mathbb{E}[w(t, X^{s,x; \alpha[r,\xi],v}_t)|\mathcal{F}_r]  +(t-r)\varphi (t-r)\ a.s.\end{equation}
Denote by $\mathcal{L}$ the set of asymptotic sub-solutions.
\end{enumerate}
\end{Definition}
\begin{Remark}\label{analytic} The definition of asymptotic solutions, tell us, that, for each $r$ there exists a Markov control at that time, such that, if the control is held constant until later, the state equation plugged inside $w$ will {\bf almost}  have  the sub-martingale property in between $r$ and any later time (reasonably close), for {\bf any choice of open loop controls} $v$. Since $v$ can change wildly after $r$, in this framework it is very convenient to consider asymptotic sub-solutions in the stochastic sense, resembling \cite{sirbu} or \cite{bs-3}. However, if the open loop controls are restricted to {\bf not} change very soon after $r$, one could consider asymptotic sub-solutions in analytic formulation. Without pursuing this direction here, in the definition of such a sub-solution, the inequality \eqref{sub-sol} could be replaced by 
$$w(r,x_1) \leq \int _{\mathbb{R}^d} w(t, z)p(r,t,x_1,x_2; \xi (x_1) ,v)dx_2 +(t-r)\varphi (t-r), \forall x_1\in \mathbb{R}^d, v\in V.$$
Here, $p(r,t,x_1,x_2,u,v)dx_2$ is the transition law from time $r$ to $t$ of the state process $X$  where $u,v$ are held constant (which is, obviously a Markov process). Such definition, as mentioned above, would work well if controls $v$ do not change in between $r$ and $t$, and would amount to "analytic asymptotic Perron's method", as opposed to the stochastic set-up we follow below.
\end{Remark}
Compared to \cite{sirbu}, the next Proposition is not entirely trivial, but not hard either.
\begin{Proposition}\label{prop1} Any $w\in \mathcal{L}$ satisfies
$w\leq v^-_M\leq v^-\leq W^-=V^-.$
\end{Proposition}
Proof: Fix $\epsilon$ and let $\delta$ such that $\varphi(\delta)\leq \varepsilon$. Choose a time partition such that $t_k-t_{k-1}\leq \delta.$ For this particular partition,  we construct, recursively, going from time $t_{k-1}$ to time $t_k$, some measurable  $\xi_k:\mathbb{R}^d\rightarrow U$ satisfying the Definition \ref{def-sub-sol}. Now, we have, with $\alpha$ formally defined as in the Definition \ref{def-sub-sol} of simple Markov strategies, that, for the simple Markov strategy $\alpha$ we have constructed, 
$$w(t_{k-1}, X^{s,x;  \alpha,v}_{t_{k-1}})\leq 
\mathbb{E}[w(t_k, X^{s,x; \alpha,v}_{t_k})|\mathcal{F}_{t_{k-1}}] +(t_k-t_{k-1})\underbrace{\varphi (t_k-t_{k-1})}_{\leq \varepsilon}\ a.s.\ \ \ \forall k.$$ 
This happens for any $x$ and any open loop control $v$. Taking expectations and summing up, we conclude that
$$w(s,x)\leq \mathbb{E}[w(T,X^{s,x;\alpha ,v}_T)] +\varepsilon \times (T-s), \forall v\in \mathcal{V}(s).$$ Taking  the infimum over $v$, since $w(T, \cdot)\leq g(\cdot)$,  we conclude that, if $|\pi|\leq \delta$ there exists 
$\alpha \in \mathcal{A}^M(s,\pi)$ such that 
$$w(s,x)\leq \inf _{v\in \mathcal{V}} \mathbb{E}[g(X^{s,x;\alpha ,v}_T)] +\varepsilon \times (T-s)\leq v^{\pi}(s,x)+\varepsilon \times (T-s) \ \  
\forall x\in \mathbb{R}^d.$$
Letting $\varepsilon \searrow 0$ we obtain the conclusion.
$\diamond$
\begin{Remark} We could let the gauge function $\varphi$ in Definition \ref{def-sub-sol} depend on the control $\xi$ (or even the time $r$) as well, to make the method even more flexible, but we just don't need that here. However, dependence on $\xi$ makes a difference if one wants to treat unbounded controls, rather than the compact case we consider here.
\end{Remark}

The next lemma is rather obvious.
\begin{Lemma} The set of asymptotic sub-solutions is directed upwards, i.e. $w_1,w_2\in \mathcal{L}$ implies
$w_1\vee w_2\in \mathcal{L}.$
\end{Lemma}
Proof: the only important thing in the proof is to notice that one can choose the gauge function $\varphi =\varphi_1\vee \varphi _2$ for $w=w_1\vee w_2$. The choice of $\xi$ is obvious. $\diamond$

\noindent {\bf Asymptotic Perron's Method for strategies:}
 we define
 $$w^-\triangleq\sup _{w\in \mathcal{L}}w\leq v^-_M\leq v^-\leq W^-.$$
 \begin{Proposition}[Asymptotic Perron]\label{perron} The function $w^-$ is an LSC viscosity super-solution of the (lower) Isaacs equation.
 \end{Proposition}
 Proof:  from Proposition 4.1 in \cite{bs-1}, there exist $\tilde{w}_n\in \mathcal{L}$ such that
$w^-=\sup _n \tilde{w}_n.$ We define the increasing sequence
$w_n =\tilde{w}_1\vee \dots \vee \tilde{w}_n\in \mathcal{L}\nearrow w^-.$

\noindent {\bf 1. Interior super-solution property}
Let $\psi$ touch $w^-$ strictly below at some $(t_0, x_0)\in [0,T)\times \mathbb{R}^d$.  Let us assume, by contradiction, that the viscosity super-solution property fails at $(t_0,x_0)$. This means that
$$\psi _t(t_0,x_0)+\sup _u \inf _v L(t_0,x_0,u,v; \psi_x (t_0,x_0), \psi _{xx}(t_0,x_0)))>0.$$
Since $L$ is continuous, we can choose a small neighborhood $B(t_0, x_0; \varepsilon) \subset [0,T)\times \mathbb{R}^d$ and some $\hat u\in U$ such that, over this neighborhood we have
$$\psi _t(t,x)+ \inf _v L(t,x,\hat{u},v; \psi_x (t,x), \psi _{xx}(t,x)))>\varepsilon.$$
From here, we follow the usual Perron construction. The very different part will be to show that, after we ``bump up'' (an approximation of) $w^-$, it still stays an {\bf asymptotic sub-solution}. More precisely, we know that, since $\psi$ touches $w^-$ below in a strict sense, there exists room of size $\delta >0$ in between $w^-$ and $\psi$ over the compact (rectangular) torus  
$$\mathbb{T}\triangleq \overline{B(t_0, x_0, \varepsilon)}-B(t_0, x_0, \varepsilon/2),$$
i.e. $w^- \geq \psi + \delta $ on $\mathbb{T}$. A Dini type argument (see, for example, \cite{bs-2}) shows that, one of the terms $w\triangleq w_n$ actually satisfies $w\geq \psi +\delta /2$ on $\mathbb{T}$. Define now, for $0<\rho <<\delta /2$ the function
$$\hat{v}=\left \{
\begin{array}{ll}
w\vee (\psi +\rho),\ \ \textrm{on}\ B(t_0, x_0; \varepsilon) \\
w,\ \textrm{outside}\ B(t_0, x_0; \varepsilon).
\end{array}\right.$$
Note that $\hat{v}=w$ on  the overlapping $\mathbb{T}$ (so, it is continuous) and $\hat{v}(t_0,x_0)=w^-(t_0,x_0)+\rho >w^-(t_0,x_0)$. The proof would be finished if we can show that $\hat{v}$ is an asymptotic sub-solution.

\noindent {\bf The idea of the proof} is quite simple, namely:
\begin{enumerate}
\item if, at time $r$, we have $w\geq \psi +\rho $ (at that particular position $y(r)$), then, we follow from $r$ forward the nearly optimal strategy corresponding to the asymptotic sub-solution $w$ (which depends only on $y(r)$)
\item if, at time $r$, we have instead $w<\psi +\rho$ (again, at that particular position $y(r)$) we follow from $r$ forward the strategy $\hat{u}$. In between $r$ and any later time $t$, the process $\psi+\gamma$ super-posed to the state equation is not a true sub-martingale, but is an asymptotic one. The reason is that, it is a sub-martingale  until the first time it exits $B(t_0,x_0;\varepsilon)$. However, the chance that this happens before $t$ can be  estimated in terms of the size of the interval $t-r$, and bounded above by a gauge function.
\end{enumerate} 
We develop rigorously below the arguments described above.
Fix $s\leq r\leq T$. Since $w$ is an asymptotic sub-solution, there exists a Markov strategy $\xi$ at time $r$ corresponding to  the Definition \ref{def-sub-sol}  for the sub-solution $w$ (for the initial time $s$). Now, we define
$$\hat{\xi}(x)=1_{\{(r,x)\notin B(t_0,x_0;  \varepsilon/2) \vee  w (r,x)\geq \psi(r,x) +\rho\}}\xi (x)+1_{\{ (r,x)\in B(t_0,x_0 ; \varepsilon/2)\wedge  w (r,x)< \psi(r,x) +\rho\}} \hat{u}.$$
We want to show that $\hat{\xi}$ satisfies the desired property in the Definition  \ref{def-sub-sol} for the (expected)  sub-solution $\hat{v}$ at $r$, with an appropriate choice of the gauge function $\varphi$ {\bf independent of} $r$,  $s$ or $\hat{\xi}$. Let $\varphi _w$ be the gauge function of the sub-solution $w$.
Consider any $\alpha \in \mathcal{A}(s)$ and any $v\in \mathcal{V}(s)$. By the definition of the sub-solution $w$, and, taking into account that
$w\leq \hat{v}$ we have that, on the event 
$$A\triangleq \{(r,X^{s,x;\alpha,v}_r)\notin B(t_0,x_0  ;\varepsilon/2)\} \cup  \{w (r,X^{s,x;\alpha, v}_r
)\geq \psi(r,X^{s,x;\alpha , v}_r) +\rho\}\in \mathcal{F}_r$$
we have that $X^{s,x;\alpha [r,\xi],v}_t=X^{s,x;\alpha [r,\hat{\xi}],v}_t$ a.s. for $r\leq t\leq T$ and, therefore
\begin{equation}
\label{1}
\begin{split}
1_A \hat{v}(r, X^{s,x; \alpha,v}_r)=   1_A w(r, X^{s,x; \alpha,v}_r)\leq &
 \mathbb{E}[1_Aw(t, X^{s,x; \alpha[r,\xi],v}_t)|\mathcal{F}_r] +\\+1_A\times (t-r)\varphi _w (t-r)
\leq& \mathbb{E}[1_A \hat{v}(t, X^{s,x; \alpha[r, \hat{\xi}],v}_t)|\mathcal{F}_r] +1_A\times (t-r)\varphi _w (t-r).
\end{split}
\end{equation}
On the complement of $A$, the process $\psi (t, X^{s,x;\alpha[r,\hat{\xi}],v}_t)$ is a sub-martingale (by It\^o) up to the first time $\tau$ where the process gets out of $B(t_0, x_0; \varepsilon)$, i.e. up to 
$$\tau\triangleq \inf \{t\geq r|  (t, X^{s,x;\alpha[r,\hat{\xi}]}_t) \notin  B(t_0, x_0; \varepsilon)\}.$$
The sub-martingale property says that
$$1_{A^c}  (\psi+\rho)(r, X^{s,x; \alpha,v}_r) \leq \mathbb{E}[1_{A^c}(\psi +\rho) (\tau \wedge t, X^{s,x; \alpha[r,\hat{\xi}],v}_{\tau \wedge t})|\mathcal{F}_r]\leq  \mathbb{E}[1_{A^c} \hat{v} (\tau \wedge t, X^{s,x; \alpha[r,\hat{\xi}],v}_{\tau \wedge t})|\mathcal{F}_r].
$$
 Fix $t$ such that $r\leq t\leq r+\varepsilon/2$. Denote now the event
$$B\triangleq \{ |X^{s,x;\alpha[r,\hat{\xi}],v}_ {t'}-x_0|<\varepsilon, \forall \ r\leq t'\leq t\}.$$
We use here the norm $|(t,x))\triangleq \max \{|t|, |x|\}.$ Consequently, we have 
$$ \mathbb{E}[1_{A^c }\hat{v} (\tau \wedge t, X^{s,x; \alpha[r,\hat{\xi}],v}_{\tau \wedge t})|\mathcal{F}_r]=
 \mathbb{E}[1_{A^c} 1_B \hat{v} (t, X^{s,x; \alpha[r,\hat{\xi}],v}_t)|\mathcal{F}_r]+
 \mathbb{E}[1_{A^c} 1_{B^c}\hat{v} (\tau, X^{s,x; \alpha[r,\hat{\xi}],v}_{\tau})|\mathcal{F}_r].
$$
Therefore, 
\begin{equation*}\begin{split}
\mathbb{E}[1_{A^c} \hat{v} (\tau \wedge t, X^{s,x; \alpha[r,\hat{\xi}],v}_{\tau \wedge t})|\mathcal{F}_r]=&
\mathbb{E}[1_{A^c} \hat{v} (t, X^{s,x; \alpha[r,\hat{\xi}],v}_t)|\mathcal{F}_r]+ \\ + &  \mathbb{E}[1_{A^c} 1_{B^c}\big(\hat{v} (\tau, X^{s,x; \alpha[r,\hat{\xi}],v}_{\tau}) -\hat{v} (t, X^{s,x; \alpha[r,\hat{\xi}],v}_{t}) \big)|\mathcal{F}_r].
\end{split}
\end{equation*}
Since $\hat{v}$ is bounded by some constant $\|\hat{v}\|_{\infty}$, we conclude that, for $r\leq t\leq r+\varepsilon/2$ we have
\begin{equation}\label{2}\begin{split} 1_{A^c}  \hat{v}(r, X^{s,x; \alpha,v}_r)=&1_{A^c}  (\psi+\gamma)(r, X^{s,x; \alpha,v}_r) \leq\\
\leq &
\mathbb{E}[1_{A^c} \hat{v} (t, X^{s,x; \alpha[r,\hat{\xi}],v}_t)|\mathcal{F}_r]+
2\|\hat{v}\|_{\infty}\mathbb{P}[A^c\cap B^c|\mathcal{F}_r], a.s. \end{split} \end{equation}
We can now put together \eqref{1} and \eqref{2}.
If we can find a gauge function $\tilde{\varphi}$ such that
$$2\|\hat{v}\|_{\infty}\mathbb{P}[A^c\cap B^c|\mathcal{F}_r]\leq 1_{A^c}\tilde{\varphi}(t-r)\times (t-r), \ a.s.$$
we are done, as one can choose the gauge function for $\hat{v}$ as
$$\varphi _{\hat{v}}\triangleq \varphi _w\vee \tilde{\varphi}.$$
We do that in the Lemma \ref{prob-estimate} below, and finish the proof of the interior sub-solution property.
\begin{Lemma}\label{prob-estimate}  There exist constants $C,C'$ (depending only on  the function $v$) such that, for any $s$ and any $r\geq s$, if $t\geq r$ is close enough to $r$ we have 
$$\mathbb{P}[A^c\cap B^c|\mathcal{F}_r]=
1_{A^c}\mathbb{P}[A^c\cap B^c|\mathcal{F}_r]  \leq 1_{A^c}C'\mathbb{P}\left (N(0,1)\geq \frac 1{C\sqrt{t-r}} \frac{\varepsilon}4 \right) \ \ a.s,$$
independently over all strategies  $\alpha \in \mathcal{A}(s)$ and controls $v\in \mathcal{V}$.
The function $$\tilde{\varphi}(t)\triangleq \frac{2\|v\|_{\infty}C'}t \mathbb{P}\left (N(0,1)\geq \frac 1{C\sqrt{t}} \frac{\varepsilon}4 \right)$$satisfies
$\lim _{t\searrow 0}\tilde{\varphi}(t)=0$ and therefore is a gauge function.

\end{Lemma}
Proof:   To begin with, we emphasize that we do not need such a precise bound on conditional probabilities, to finish the proof of Theorem \ref{main} (both the interior super-solution part or the terminal condition). The simple idea of the proof is to see that, conditioned on $A^c$, the event we care about amounts to a continuous semi-martingale with bounded volatility  and bounded drift to exit from a {\bf fixed} box in the interval  of time $[r,t]$. If the size of $t-r$ is small enough, that amounts to just the martingale part exiting from a smaller fixed box in between $t$ and $r$. This can be rephrased, through a time change, in terms of a Brownian motion, and estimated very precisely to be of the order $ \mathbb{P}\left (N(0,1)\geq \frac 1{C\sqrt{t-r}} \frac{\varepsilon}4 \right)=o(t-r)$, where $N(0, a^2)$ is a normal with mean zero and standard deviation $a$.  This is basically the whole proof, in words. The  precise mathematics below follows exactly these lines.
We first notice that
$$A^c\cap B^c\subset {\{(r,X^{s,x;\alpha, v}_r)\in B(t_0,x_0  ;\varepsilon/2)   \ \textrm{and}\  X^{s,x;\alpha[r,\hat{\xi}],v}_ {t'}\notin B(t_0, x_0; \varepsilon),  \ \textrm{ for \ some}\  r\leq t'\leq t\}}.$$ If $t-r< \varepsilon/2$ we  have
$$A^c\cap B^c\subset {\{(r,X^{s,x;\alpha, v}_r)\in B(t_0,x_0  ;\varepsilon/2)   \ \textrm{and}\  |X^{s,x;\alpha[r,\hat{\xi}],v}_ {t'}-  X^{s,x;\alpha [r, \xi]}_r|\geq \varepsilon/2
,  \ \textrm{ for \ some}\  r\leq t'\leq t\}}.$$
Over $B(t_0, x_0;\varepsilon)$ both the drift and the volatility of the state system are uniformly bounded by some constant $C$. Therefore, if we choose 
$t-r\leq \frac{\varepsilon}{4C}$ the integral of the drift part cannot exceed $\varepsilon /4$ in size. We, therefore, conclude that,with the notation
$$D\triangleq \left \{ (r,X^{s,x;\alpha,v}_r)\in B(t_0,x_0  ;\varepsilon/2)\right\},$$
if  $t-r\leq \frac{\varepsilon}2\wedge \frac{\varepsilon}{4C},$ then
$$A^c\cap B^c\subset D\cap \left \{  \left |  \int _r^{t'} \sigma (q, X^{s,x;\alpha[r, \hat{\xi}],v}_q,\alpha[r, \hat{\xi}](q, X^{s,x;\alpha[r, \hat{\xi}],v}_{\cdot}), v_q)dW_q\right| \geq \varepsilon/4 \ \textrm{  for \  some}\  r\leq t'\leq t\right\}.$$
Denote by 
$$M_{t'}\triangleq  \int _s^{t'} \sigma (q, X^{s,x;\alpha[r, \hat{\xi}],v}_q,\alpha[r, \hat{\xi}](q, X^{s,x;\alpha[r, \hat{\xi}],v}_{\cdot}), v_q )dW_q, \forall  \ \ s\leq t'\leq T.$$
We study separately the coordinates. More precisely 
we consider on $\mathbb{R}^d$ the max norm as well, and, for 
$M_{t'}=(M^1_{t'},\dots, M^d_{t'})$, we also have
$$A^c\cap B^c\subset \bigcup _{l=1}^d \left (D\cap \{ \rho_l\leq t\}\right), $$
where
$$\rho_l \triangleq \inf \{r\leq t' \leq \tau  | |M^l_{t'}-M^l_r|\geq \varepsilon/4\}.$$
We estimate the probabilities above (conditioned on $\mathcal{F}_r$) individually.
We want to use  the result of Dambis-Dubins-Schwarz to translate the computation into a probability of exiting from a fixed box of a Brownian Motion. In order to do so rigorously, and without enlarging the probability space to accommodate an additional Brownian motion, we first need to make the volatility  explode at $T$. Choose any function 
 $$f\rightarrow [s,T)\rightarrow (0,\infty),\ \ \int _s^tf^2(q)dq<\infty \ \forall q<T, \   \int_s^Tf^2 (q)dq=\infty. $$ 
 Fix some $\varepsilon '>0$ such that $t_0+\varepsilon <T-\varepsilon '$, so $\tau \leq T-\varepsilon'$ on $D$. 
 Define
 $$M^{\varepsilon'}_{t'}\triangleq  \int _s^{t'} \left (1_{\{s\leq q\leq T-\varepsilon'\}}\sigma (q, X^{s,x;\alpha[r, \hat{\xi}],v}_q,\alpha[r, \hat{\xi}](q, X^{s,x;\alpha[r, \hat{\xi}],v}_{\cdot}), v_q )+1_{\{T-\varepsilon '<q\leq T\}}f(q) \right )dW_q, \forall s\leq t' <T.$$ 
 We have
$$\rho_l =\inf \{r\leq t' \leq \tau  | |M^{\varepsilon',l}_{t'}-M^{\varepsilon',l}_r|\geq \varepsilon/4\}.$$
 According to Dumbis-Dubins-Schwarz (see \cite{KS88}, page 174) applied with the obvious shift of time origin from $r$ to zero, 
 the process 
  $(Z_{s'})_{0\leq s'<\infty}$ defined as
 $$Z_{s'}=M^{\varepsilon',l}_{A_{s'}}-M^{\varepsilon',l}_r, 0\leq s'<\infty,\ \textrm{ for}\ 
A_{s'}=\inf \{t'\geq r|\langle M^{\varepsilon',l}\rangle _{t'} -\langle M^{\varepsilon',l}\rangle _r>s'\}.$$
 is a Brownian motion with respect to the filtration $\mathcal{G}_{s'}\triangleq \mathcal{F}_{A_{s'}}, \ 0\leq s'<\infty.$
Therefore, since $\mathcal{F}_r\subset \mathcal{G}_0$, $Z$  is 
  {\bf independent} of $\mathcal{F}_r$. In addition, still from DDS, we have
 $$M^{\varepsilon',l}_{t'}-M^{\varepsilon',l}_r=Z_{\langle M^{\varepsilon',l}\rangle_{t'}-\langle M^{\varepsilon',l}\rangle _r}, r\leq t'< T.$$ 
 On the other hand,  since $\sigma$ is uniformly bounded (independently on the strategy and the control) on $D$ and before $\tau$ (the exit time from $B(t_0,x_0;\varepsilon)$)  occurs we have 
 $$\langle M^{\varepsilon',l}\rangle_{t'}-\langle M^{\varepsilon',l}\rangle _r \leq C^2(t'-r),\ \ r\leq t'\leq \tau \leq T-\varepsilon '.$$ 
 We then conclude that
 $$D\cap \{ \rho_l\leq t\}  
 \subset \left\{ \max _{0\leq s'\leq C^2(t-r) }|Z_{s'}|\geq \varepsilon/4\right\}.$$
 Let $M^+_t$ and $M^-_t$ the distributions of the running max and running min of a standard BM starting at time zero.
 Since the Brownian motion $Z$ is independent of $\mathcal{F}_r$, we can obtain that
 $$\mathbb{P}(D\cap \{ \rho_l\leq t\}  |\mathcal{F}_r)\leq \mathbb{P}(M^+_{C^2(t-r)}\geq \varepsilon/4 ) + \mathbb{P}(M^-_{C^2(t-r)}\geq \varepsilon/4 ) =2\mathbb{P}(M^+_{C^2(t-r)}\geq \varepsilon/4 ).$$
 It is well know, from the reflection principle, (see, for example, Karatzas and Shreve \cite{KS88} page 80) that this is the same as
 $$\mathbb{P}(D\cap \{ \rho_l\leq t\}  |\mathcal{F}_r)\leq4\mathbb{P}(N(0,C^2(t-r))\geq \varepsilon /4)=4\mathbb{P}\left (N(0,1)\geq \frac 1{C\sqrt{t-r}} \frac{\varepsilon}4 \right).$$
 We emphasize that all estimates are {\bf independent} of the initial time $s$, the later time $r$ and the even later time $t$  as long as  $t-r\leq \frac{\varepsilon}2\wedge \frac{\varepsilon}{4C}$  (and  independent of $\alpha \in \mathcal{A}(s)$ and of $v\in \mathcal{V}(s)$).  
Summing all the terms, we obtain
 $$\mathbb{P}[A^c\cap B^c|\mathcal{F}_r]\leq  4d \mathbb{P}\left (N(0,1)\geq \frac 1{C\sqrt{t-r}} \frac{\varepsilon}4 \right), \ \textrm{if}\  t-r\leq \frac{\varepsilon}2\wedge \frac{\varepsilon}{4C}.$$  We can multiply now with $1_{A^c}$ which is measurable with respect to $\mathcal{F}_r$ to obtain the conclusion. It is well known that
 $$\frac{ \mathbb{P}\left (N(0,1)\geq \frac 1{C\sqrt{t}} \frac{\varepsilon}4 \right)}{t}\rightarrow 0,\ \ \textrm{as}\ \ t\searrow 0,$$
 and this finishes the proof of the lemma. Instead of appealing to Dambis-Dubins-Scwhwarz theorem, we could also use the (smooth solution) characterizing the exit probably of a Brownian motion from a box,  and super-pose it to the state process, resulting in a super-martingale (just using It\^o formula). This would still bound the probability we are interested in by the exit probability of a standard Brownian motion from a box.$\diamond$

\noindent Proof of Proposition \ref{perron}, continued: 

 \noindent {\bf 2. The terminal condition $w^-(T,\cdot)\geq g(\cdot)$} The proof of this is done again, by contradiction. The "bump-up" analytic construction we use is similar to \cite{bs-3}, and the rest is based on similar arguments to the interior super-solution property and a very  similar estimate to Lemma \ref{prob-estimate} above.  The only difference, if the Dambis-Dubins-Schwarz route is followed, is to  see that, with all notations as above, if $r\leq t<T,$ then 
 $$\{\rho_l\leq t\} =\lim _{\varepsilon'\searrow 0} \{\rho_l \leq t\wedge (T-\varepsilon')\}=
 \lim _{\varepsilon'\searrow 0} \{|M^{\varepsilon',l}_{t'}-M^{\varepsilon',l}_r|\geq \varepsilon/4\ \ \textrm{for  \ some}\ r\leq t'\leq \tau \wedge t\wedge(T-\varepsilon')\} .$$
 We can first let $\varepsilon '\searrow 0$ to obtain the conclusion that (with the same notation $\hat{v}$ for the "bump-up" function), with the gauge function we just chose in the proof of the Lemma above we have
 $$
 \hat{v}(r, X^{s,x; \alpha,v}_r) \leq 
\mathbb{E}[\hat{v}(t, X^{s,x; \alpha[r,\hat{\xi}],v}_t)|\mathcal{F}_r]  +(t-r)\varphi (t-r)\ a.s$$ 
for all $t\in [r,T)$ close enough to $r$, i.e. $t-r\leq \frac{\varepsilon}2\wedge\frac{\varepsilon}{4C}.$ 
Obviously, letting $t\nearrow T$ we obtain the same for $t=T$. The proof is now complete.
 $\diamond$ 
 
 \subsection{Asymptotic Perron over counter-strategies}\label{subsection-counter}
One only has to go through this construction if the genuine  (non-symmetric) game  in \eqref{non-sym-game} is studied. The notion of counter-strategies (even with the Markov discretization) is still not so easily implementable, since actions do change continuously in time (if  $u$ does so, in a situation of  counter-strategies  vs. open loop-controls as in the definition \eqref{v+} of of $v^+$, which describes some strange model of  "worst case scenario", analyzed here for mathematical  reasons only). On the other  hand, there is no way one can genuinely discretize counter-strategies to obtain a value in \eqref{non-sym-game} (see Remark \ref{counter-discrete}). Therefore, we go over this analysis for mathematical completeness, emphasizing that the basic  method and the more important result are contained in the previous Subsection \ref{subsection-strategies}.
We view this as a  simple additional application of the Asymptotic Perron  method. The definitions and  proofs follow in lockstep with the previous Subsection \ref{subsection-strategies}, with minor appropriate modifications to account for counter-strategies.

\begin{Definition}[Asymptotic Stochastic Super-Solutions]\label{def-super-sol}
A function $w:[0,T]\times \mathbb{R}^d\rightarrow \mathbb{R}$ is called an asymptotic (stochastic) super-solution of the (lower) Isaacs equation, if it is bounded, continuous and  satisfies $w(T,\cdot)\geq g(\cdot)$. In addition, there exists a gauge function $\varphi=\varphi_w:(0,\infty)\rightarrow (0,\infty),$ depending on $w$ such that
\begin{enumerate}
\item $\lim _{\varepsilon\searrow 0}\varphi(\varepsilon)=0,$
\item for each $s$ (and the optimization problem coming with it),  for each time $s\leq r\leq T$, there exists a measurable function $\eta :\mathbb{R}^d  \times U \rightarrow V$ such that, for each $x$, each  $\gamma \in \mathcal{C}(s)$ and $u\in \mathcal{U}(s)$, if we make the notation
$\gamma [r,\eta]\in \mathcal{C}(s),$
defined by
$$\gamma [r,\eta](t, y(\cdot),u)=1_{\{s<t\leq r\}}\gamma (t, y(\cdot), u)+1_{\{r<t\leq T\}} \eta (y(r), u),$$
then, for each $r\leq t\leq T$ we have
\begin{equation}\label{sub-sol}
w(r, X^{s,x; u,\gamma}_r)=w(r, X^{s,x; u, \gamma [r,\eta]}_r)\geq 
\mathbb{E}[w(t, X^{s,x; u, \gamma [r,\eta],}_t)|\mathcal{F}_r]  -(t-r)\varphi (t-r)\ a.s.\end{equation}
Denote by $\mathcal{U}$ the set of asymptotic super-solutions.
\end{enumerate}
\end{Definition}
We again have:
\begin{Proposition}\label{prop2} Any $w\in \mathcal{U}$ satisfies
$w\geq v_M^+$.
\end{Proposition}
Proof: Fix $\epsilon$ and let $\delta$ such that $\varphi(\delta)\leq \varepsilon$. Choose $\pi$  such that $\|\pi\|\leq \delta.$ For this partition $\pi$,  we construct, recursively, going from time $t_{k-1}$ to time $t_k$, some measurable  $\eta_k:\mathbb{R}^d \times U \rightarrow V$ satisfying the Definition \ref{def-super-sol}. 
 We put the $\eta_k$'s together to obtain a Markov counter-strategy $\gamma$ for which
$$w(t_{k-1}, X^{s,x;  u, \gamma}_{t_{k-1}})\geq 
\mathbb{E}[w(t_k, X^{s,x; u, \gamma }_{t_k})|\mathcal{F}_{t_{k-1}}] -(t_k-t_{k-1})\underbrace{\varphi (t_k-t_{k-1})}_{\leq \varepsilon}\ a.s.\ \ \ \forall k.$$ 
This happens for any $x$ and any open loop control $u$. Taking expectations and summing up, we conclude that
$$w(s,x)\geq \mathbb{E}[w(T,X^{s,x;u ,\gamma}_T)] -\varepsilon \times (T-s), \forall u\in \mathcal{U}(s).$$ Taking  the supremum over $u$, since $w(T, \cdot)\geq g(\cdot)$,  we conclude that, if $|\pi|\leq \delta$ there exists 
$\gamma \in \mathcal{C}^M(s,\pi)$ such that 
$$w(s,x)\geq \sup _{u\in \mathcal{U}(s)} \mathbb{E}[g(X^{s,x;u ,\gamma}_T)] -\varepsilon \times (T-s)\geq v_{\pi}^+(s,x)-\varepsilon \times (T-s) \ \  
\forall x\in \mathbb{R}^d.$$
Letting $\varepsilon \searrow 0$ we obtain the conclusion.
$\diamond$

The next lemma is, once again,  obvious.
\begin{Lemma} The set of asymptotic super-solutions is directed downwards, i.e. $w_1,w_2\in \mathcal{U}$ implies
$w_1\wedge  w_2\in \mathcal{U}.$
\end{Lemma}
Proof: the only thing needed is  to notice that one can choose the gauge function $\varphi =\varphi_1\vee \varphi _2$ for $w=w_1\vee w_2$. The choice of $\eta$ is, again,  obvious. $\diamond$

\noindent {\bf Asymptotic Perron's Method for counter-strategies:}
 we define
 $$w^+\triangleq\inf  _{w\in \mathcal{U}}w\geq v_+^M \geq v^+\geq W^+.$$
 \begin{Proposition}[Asymptotic Perron]\label{perron-counter} The function $w^+$ is an USC viscosity sub-solution of the (lower) Isaacs equation.
 \end{Proposition}
 Proof:
\noindent {\bf 1. Interior sub-solution property}
Let $\psi$ touch $w^+$ strictly above at some $(t_0, x_0)\in [0,T)\times \mathbb{R}^d$.  Assume, by contradiction, that 
$$\psi _t(t_0,x_0)+\sup _u \inf _v L(t_0,x_0,u,v; \psi_x (t_0,x_0), \psi _{xx}(t_0,x_0)))<0.$$
This means that there exists a  small $\varepsilon >0$ and a (measurable)
 function 
$h: U\rightarrow V$  such that
$$\psi _t(t_0,x_0)+ L(t_0,x_0,u,h(u); \psi_x (t_0,x_0), \psi _{xx}(t_0,x_0)))<-\varepsilon.$$
Since $L$ is continuous (so uniformly continuous over $(t,x,u,v,p,M)$ as long as  $(t,x)$ is close to $(t_0, x_0)$ and $(p.M)$ is close to $(\psi _x(t_0, x_0), \psi _{xx}(t_0, x_0))$), we can choose an even smaller $\varepsilon$ 
such that 
$$\psi _t(t,x)+  L(t,x,u,h(u); \psi_x (t,x), \psi _{xx}(t,x)))<-\varepsilon$$
over the (smaller) neighborhood $B(t_0, x_0; \varepsilon) \subset [0,T)\times \mathbb{R}^d$. From here on, we  follow the usual Perron construction. We need to show that after we ``bump down'' (an approximation of) $w^+$, it still stays an {\bf asymptotic super-solution}. 
Since $\psi$ touches $w^+$ above in a strict sense, there exists room of size $\delta >0$ in between $w^+$ and $\psi $ over the compact (rectangular) torus  
$$\mathbb{T}\triangleq \overline{B(t_0, x_0, \varepsilon)}-B(t_0, x_0, \varepsilon/2),$$
i.e. $w^+ \leq \psi - \delta $ on $\mathbb{T}$. Again a  Dini type argument (see, for example, \cite{bs-2}) shows that, one of the terms  of the sequence $w_n \searrow w^+$, which we simply denote by $w$, actually satisfies $w\leq \psi -\delta /2$ on $\mathbb{T}$. Define now, for $0<\rho <<\delta /2$ the function
$$\hat{v}=\left \{
\begin{array}{ll}
w\wedge  (\psi -\rho),\ \ \textrm{on}\ B(t_0, x_0; \varepsilon) \\
w,\ \textrm{outside}\ B(t_0, x_0; \varepsilon).
\end{array}\right.$$
We have that $\hat{v}=w$ on  the overlapping $\mathbb{T}$ (so,  $\hat{v}$ it is continuous) and $\hat{v}(t_0,x_0)=w^+(t_0,x_0)-\rho <w^+(t_0,x_0)$. We only need to show that $\hat{v}$ is an asymptotic super-solution to have a full proof. Fix $s\leq r\leq T$. Since $w$ is an asymptotic super-solution in the sense of Definition \ref{def-super-sol}, there exists an $\eta :\mathbb{R}^d\times U\rightarrow V$ at time $r$ corresponding to  the Definition \ref{def-super-sol}  for the super-solution $w$ (for the initial time $s$). Define
\begin{equation}
\label{counter-super-opt}\hat{\eta}(x,u)=1_{\{(r,x)\notin B(t_0,x_0;  \varepsilon/2) \vee  w (r,x)\leq \psi(r,x) -\rho\}}\eta (x,r)+1_{\{ (r,x)\in B(t_0,x_0 ; \varepsilon/2)\wedge  w (r,x)> \psi(r,x) -\rho\}} h(u).
\end{equation}
\begin{Remark}
The choice of $h(u)$ together with It\^o formula tells as that, as long as the player $v$ always adjusts his/her control (observing continuously the other player's actions) to be $h(u_t)$ and  $(t,X_t)$ is inside $ B(t_0, x_0; \varepsilon)$ then 
$(\psi -\rho )(t,X_t)$ is a super-martingale.\end{Remark}
Now, some very similar arguments  to the considerations in Subsection \ref{subsection-strategies}  based on the Remark above and a  (next to) identical result to Lemma \ref{prob-estimate} finish the proof of $\hat{v}\in \mathcal{U}$, resulting in a contradiction. $\diamond$

\noindent{\bf 2. The Terminal condition $w^+(T,\cdot)\leq g(\cdot)$}
One has to use identical arguments to the proof of  $w^-  (T,\cdot)\geq g(\cdot)$ in Subsection \ref{subsection-strategies} (which was, in turn, apparent), with the only difference of constructing a counter-strategy similar to \eqref{counter-super-opt} in order to reach a contradiction to 
the assumption that  $w^+(T,x_0)> g(x_0)$ for some $x_0$.

 \subsection{Proof of Theorem \ref{main}} 
 Recall that
 $w^-\leq v^-_M\leq v^-\leq W^- =V^-\leq W^+\leq v^+\leq v^+_M\leq w^+.$
  Since $w^-$ is a LSC viscosity super-solution and $w^+$ is an USC sub-solution of the lower Isaacs equation, 
 the comparison result in \cite{sirbu} ensures that 
$$v^-_M=v^-= W^-=V^-= W^+= v^+=v^+_M.$$
  As mentioned before, if one does not really care about the genuine non-symmetric game \eqref{non-sym-game} and it's value/saddle points, then  only the  the lower Perron construction 
$$w^-\leq v^-_M\leq v^-\leq V^-$$  is needed.
The viscosity super-solution property of $w^-$ together with the viscosity property of $V^-$ from \cite{sirbu} (which is actually re-proved above) yields  the more important half of the Theorem \ref{main}, which is 
  $$ v^-_M= v^-= V^-.$$  
Now, in order to prove the second part of Theorem \ref{main}, we note that we have constructed $\mathcal{L}\ni w^n\nearrow W^-$. By continuity and the Dini's criterion, the above convergence is uniform on compacts.  This means that, for each $\varepsilon$ there exists $w\in \mathcal{L}$ (one of the terms of the increasing sequence of asymptotic sub-solutions) such that
$$ W^- -\varepsilon \leq w, \ \textrm{on}\ C=[0,T]\times \ \{|x|\leq N\}.$$
Let $\varphi$ be the gauge function of this particular $w$, and let $\delta$ such that $\varphi (\delta)\leq \varepsilon.$
According to the proof of Proposition \ref{prop1}, if $|\pi|\leq \delta$, there exists
$\hat{\alpha }\in \mathcal{A}^M(s,\pi)$ such that
$$w(s,x)\leq \inf _{v\in \mathcal{V}} \mathbb{E}[g(X^{s,x;\hat{\alpha} ,v}_T)] + \varepsilon \times (T-s)\ \forall x.$$ This implies that
$$ W^-(s,x) -\varepsilon\times (1+(T-s)) \leq \inf _{v\in \mathcal{V}(s)} \mathbb{E}[g(X^{s,x;\hat{\alpha },v}_T)] \leq v^- _{\pi} (s,x),\ \ \forall |x|\leq N.$$
A very similar argument based on Dini, together with Proposition \ref{prop2} shows that, for $\|\pi\|\leq \delta (\varepsilon)$ (here $\delta (\varepsilon)$ may have to be modified) there exists a counter-strategy $\hat{\gamma }\in \mathcal{C}(\pi,s)$ such that 
$$ W^+(s,x) +\varepsilon\times (1+(T-s)) \geq \sup _{u\in \mathcal{U}(s)} \mathbb{E}[g(X^{s,x;u,\hat{\gamma}}_T)] \geq v^+_{\pi}(s,x),\ \ \forall |x|\leq N.$$
Not only that the approximations are uniform on $C$, but, for fixed time $s$, the uniform approximations can be realized over the same simple Markov strategy $\hat{\alpha} \in \mathcal{A}^M(s,\pi)$ or the same Markov counter-strategy $\hat{\gamma}\in \mathcal{C}^M(s,\pi)$ for $|\pi |\leq \delta.$
$\diamond$
\section{Final Considerations}
\subsection{One player/control problems}\label{one-player}
In case the state system only depends on $u$ and not on $v$ (i.e. we have a control problem rather than a game), then, with the obvious observation that
$$v_M^-(s,x)\leq V ^-(s,x)\leq \sup_{u\in \mathcal{U}(s)} \mathbb{E}[g(X^{s,x;u}_T)], $$
one can use our result about games to conclude that, in a control problem (one-player) like in \cite{bs-3} (but under our stronger standing assumptions here), the value functions over open-loop controls, elementary feed-back strategies and simple Markov strategies  coincide. In addition, the approximation with simple Markov strategies is uniform over the mesh of the grid, uniform on compacts. We remind the reader that, in \cite{bs-3}, the value function studied was defined over open-loop controls, i.e.
$$V_{ol}(s,x)\triangleq \sup_{u\in \mathcal{U}(s)}{E}[g(X^{s,x;u}_T)]. $$ In this case, as pointed out in Remark \ref{analytic}, we can actually use the analytic formulation of asymptotic solutions.
Up to some considerations related to the Markov property of SDE's and some other small technical considerations (filtration, and local Lipschitz condition), this result is the same as Theorem 2 on page 148 in the seminal monograph \cite{MR2723141}. Again, we just present a novel method to prove such a result.

\subsection{Values for symmetric  feedback games }
In the case of symmetric feedback games,  if the Isaacs condition is satisfied, we know from \cite{sirbu} that the game has a value. Applying the asymptotic Perron method over strategies in  Subsection \ref{subsection-strategies} to both players (on both sides) we obtain that, for each $\varepsilon$, there exist $\varepsilon$-saddle point within the class of simple Markov strategies, uniformly in bounded $x$, which means
$(\alpha (\varepsilon),\beta(\varepsilon))\in \mathcal{A}^M(s)\times \mathcal{B}^M(s)$
such that 
$$\mathbb{E}[g(X^{s,x;u, \beta (\varepsilon)}_T]-\varepsilon\leq \mathbb{E}[g(X^{s,x;\alpha (\varepsilon), \beta (\varepsilon)}_T]
 \leq 
 \mathbb{E}[g(X^{s,x;\alpha (\varepsilon), v }_T]
  +\varepsilon\  (\forall) \ \ (u,v) \in \mathcal{U}(s)\times \mathcal{V}(s), |x|\leq N.$$
  
  If the Isaacs condition fails, we can still model the game, in a martingale formulation, as in \cite{sirbu-2}, and a value over feed-back mixed/relaxed strategies does exist. Using again the Asymptotic Perron's method, for both players, we can obtain the existence of $\varepsilon$-saddle point  within the class of mixed/relaxed strategies of simple Markov type, uniformly in bounded $x$. A mixed strategy $\mu$  of simple Markov type (for the player $u$) is defined by 
 a time grid $\pi$  and 
 some functions $\xi_k:\mathbb{R}^d\rightarrow \mathcal{P}(U), k=1,\dots, n$ measurable, such that
$$\mu (t, y(\cdot))=\sum _{k=1}^n 1_{\{t_{k-1}<t\leq t_k\}} \xi _k(y(t_{k-1}) \in \mathcal{P}(U).$$
In other words, player $u$ decides at time $t_{k-1}$ based {\bf only} on the position at that time, what distribution he/she will be sampling {\bf continuously} from until $t_k$. Obviously, one can define similarly mixed strategies of Markov type for the $v$-player. In order to do the  analysis and obtain the approximate  mixed Markov saddle strategies, one would have to go inside the short proofs in \cite{sirbu-2} and apply Asymptotic Perron's Method for the auxiliary (and strongly defined) games in the proofs there. In other words, the above paragraph for games over pure strategies satisfying Isaacs condition applies to the auxiliary game in \cite{sirbu-2}, leading to $\varepsilon$-saddle points in the class of mixed strategies of Markov type for the original game.
  
\bibliographystyle{amsalpha} 

\begin{thebibliography}{FHH12}

\bibitem[BLQ14]{buc-li-q}
R.~Buckdahn, J.~Li, and M.~Quincampoix, \emph{Value in mixed strategies for
  zero-sum stochastic differential games without {I}saacs condition}, Annals of
  Probability \textbf{42} (2014), no.~4, 1724--1768.

\bibitem[BN]{bouchard-nutz}
B.~Bouchard and M.~Nutz, \emph{Stochastic target games and dynamic programming
  via regularized viscosity solutions}, http://arxiv.org/abs/1307.5606.

\bibitem[BS91]{barles-souganidis}
G.~Barles and P.E. Souganidis, \emph{Convergence of approximation schemes for
  fully nonlinear second order equation}, Asymptotic Analysis \textbf{4}
  (1991), 271--283.

\bibitem[BS12]{bs-1}
E.~Bayraktar and M.~S\^{i}rbu, \emph{Stochastic {P}erron's method and
  verification without smoothness using viscosity comparison: the linear case},
  Proceedings of the American Mathematical Society \textbf{140} (2012),
  3645--3654.

\bibitem[BS13]{bs-3}
\bysame, \emph{Stochastic {P}erron's method for {H}amilton-{J}acobi-{B}ellman
  equations}, SIAM Journal on Control and Optimization \textbf{51} (2013),
  no.~6, 4274--4294.

\bibitem[BS14]{bs-2}
\bysame, \emph{Stochastic {P}erron's method and verification without smoothness
  using viscosity comparison: obstacle problems and {D}ynkin games},
  Proceedings of the American Mathematical Society \textbf{142} (2014), no.~4,
  1399--1412.

\bibitem[EK72]{ek}
R.~J. Elliott and N.~J. Kalton, \emph{Values in differential games}, Bull.
  Amer. Math. Soc \textbf{72} (1972), no.~3, 427--431.

\bibitem[FHH11]{fleming-hernandez2}
Wendell~H. Fleming and Daniel Hern{\'a}ndez-Hern{\'a}ndez, \emph{On the value
  of stochastic differential games}, Commun. Stoch. Anal. \textbf{5} (2011),
  no.~2, 341--351. \MR{2814482 (2012h:91034)}

\bibitem[FHH12]{fleming-hernandez2-2}
\bysame, \emph{Strategies for differential games}, Stochastic processes,
  finance and control, Adv. Stat. Probab. Actuar. Sci., vol.~1, World Sci.
  Publ., Hackensack, NJ, 2012, pp.~89--104. \MR{2985435}

\bibitem[FS89]{fs}
W.~H. Fleming and P.~E. Souganidis, \emph{On the existence of value functions
  of two-player, zero-sum stochastic differential games}, Indiana University
  Mathematics Journal \textbf{38} (1989), no.~2, 293--314.

\bibitem[Ish87]{ishii}
H.~Ishii, \emph{Perron's method for {Hamilton-Jacobi} equations}, Duke
  Mathematical Journal \textbf{55} (1987), no.~2, 369--384.

\bibitem[JS12]{JS}
K.~Jane\v{c}ek and M.~S{\^\i}rbu, \emph{Optimal investment with high-watermark
  performance fee}, SIAM Journal on Control and Optimization \textbf{50}
  (2012), no.~2, 780--819.

\bibitem[Kry00]{krylov-shaken-coeff}
N.V. Krylov, \emph{On the rate of convergence offinite-difference
  approximations for {B}ellman's equations with variable coefficients}, Probab.
  Theory Relat. Fields \textbf{117} (2000), no.~1, 1--16.

\bibitem[Kry09]{MR2723141}
N.~V. Krylov, \emph{Controlled diffusion processes}, Stochastic Modelling and
  Applied Probability, vol.~14, Springer-Verlag, Berlin, 2009, Translated from
  the 1977 Russian original by A. B. Aries, Reprint of the 1980 edition.

\bibitem[KS88a]{KS88}
I.~Karatzas and S.~Shreve, \emph{Brownian motion and stochastic calculus},
  Springer New York, 1988.

\bibitem[KS88b]{krasovskii-subbotin-88}
N.~N. Krasovski{\u\i} and A.~I. Subbotin, \emph{Game-theoretical control
  problems}, Springer Series in Soviet Mathematics, Springer-Verlag, New York,
  1988, Translated from the Russian by Samuel Kotz. \MR{918771 (89b:90248)}

\bibitem[PZ14]{pz-game}
T.~Pham and J.~Zhang, \emph{Two person zero-sum game in weak formulation and
  path dependent {B}ellman-{I}saacs equation}, SIAM J. Control Optim
  \textbf{52} (2014), no.~4, 2090---2121.

\bibitem[S\^14a]{sirbu-3}
M.~S\^{\i}rbu, \emph{A note on the strong formulation of stochastic control
  problems with model uncertainty}, Electronic Communications in Probability
  \textbf{19} (2014), no.~81, 1--10.

\bibitem[S\^14b]{sirbu-2}
\bysame, \emph{On martingale problems with continuous-time mixing and values of
  zero-sum games without {I}saacs conditions}, SIAM Journal on Control and
  Optimization \textbf{52} (2014), no.~5, 2877--2890.

\bibitem[S\^14c]{sirbu}
\bysame, \emph{Stochastic {P}erron's method and elementary strategies for
  zero-sum differential games}, SIAM Journal on Control and Optimization
  \textbf{52} (2014), no.~3, 1693--1711.

\bibitem[{\'S}wi96a]{swiech-1}
A.~{\'S}wi{\c{e}}ch, \emph{Sub- and superoptimality principles of dynamic
  programming revisited}, Nonlinear Anal. (1996), no.~8, 1429--1436.
  \MR{1377672 (97b:49026)}

\bibitem[{\'S}wi96b]{swiech-2}
Andrzej {\'S}wi{\c{e}}ch, \emph{Another approach to the existence of value
  functions of stochastic differential games}, J. Math. Anal. Appl.
  \textbf{204} (1996), no.~3, 884--897. \MR{1422779 (97j:90091)}

\end{thebibliography}
\providecommand{\bysame}{\leavevmode\hbox to3em{\hrulefill}\thinspace}
\providecommand{\MR}{\relax\ifhmode\unskip\space\fi MR }
\providecommand{\MRhref}[2]{%
  \href{http://www.ams.org/mathscinet-getitem?mr=#1}{#2}
}
\providecommand{\href}[2]{#2}

\end{document}